\documentclass{amsart}
\usepackage{graphicx}
\newtheorem{thm}{Theorem}[section]
\newtheorem{cor}[thm]{Corollary}
\newtheorem{lem}[thm]{Lemma}
\newtheorem{prop}[thm]{Proposition}
\theoremstyle{definition}
\newtheorem{defn}[thm]{Definition}
\theoremstyle{remark}
\newtheorem{rem}[thm]{Remark}
\numberwithin{equation}{section}

\begin{document}

\title[Expander Graphs and Cayley Graphs]{On Some Expander Graphs and Algebraic Cayley Graphs}%
\author{Xiwang Cao}%
\address{Xiwang Cao is with School of Mathematical Sciences, Nanjing University of
Aeronautics and Astronautics, Nanjing 210016, email: {\tt
xwcao@nuaa.edu.cn}}%


\begin{abstract}
Expander graphs have many interesting applications in communication networks and other areas, and thus these graphs have been extensively studied in theoretic computer sciences and in applied mathematics. In this paper, we use reversible difference sets and generalized difference sets to construct more expander graphs, some of them are Ramanujan graphs. Three classes of elementary constructions of infinite families of Ramanujan graphs are provided. It is proved that for every even integer $k>4$, if $2(k+2)=rs$ for two even numbers $r$ and $s$ with $s\geq 4$ and $2s>r\geq s$, or $r\geq 4$ and $2r>s\geq r$, then there exists an $k$-regular Ramanujan graph. As a consequence, there exists an $k$-regular Ramanujan graph with $k=2t^2-2$ for every integer $t>2$. It is also proved that for every odd integer $m$, there is an $(2^{2m-2}+2^{m-1})$-regular Ramanujan graph. These results partially solved the long hanging open question for the existence of $k$-regular Ramanujan graphs for every positive integer $k$.
\end{abstract}
\maketitle
\section{Introduction}

A graph $\mathcal{X}=(V,E)$ consists of a vertex set $V$ with $|V|=n$,
an edge set $E$, and a relation that associates with each edge an pair of
vertices. A graph is finite if its vertex set and edge set are finite.
The adjacency matrix $A(\mathcal{X})$ of a graph $\mathcal{X}$ is the matrix with rows and columns
indexed by the vertices of $\mathcal{X}$ with the $uv$-entry equal to the number of edges from
vertex $u$ to vertex $v$. Notice that when $\mathcal{X}$ is
undirected, $A$ is symmetric which implies the eigenvalues of $\mathcal{X}$ are real.
The spectrum of a graph is the (multi)set of its eigenvalues. We denote by $\lambda_1(\mathcal{X}) \geq \lambda_2(\mathcal{X})\geq \cdots \lambda_n(\mathcal{X})$ the eigenvalues of
an undirected graph $A(\mathcal{X})$ and we also call them the eigenvalues of the graph $\mathcal{X }$.
Spectral graph theory is the study of the eigenvalues of graphs and it has a long history. In the early days, matrix theory and
linear algebra were used to analyze adjacency matrices of graphs. Algebraic methods
have proven to be especially effective in treating graphs which are regular and
symmetric. Sometimes, certain eigenvalues have been referred to as the ``algebraic
connectivity" of a graph. There is a large literature on algebraic aspects of
spectral graph theory, such as Biggs
\cite{biggs}, Cvetkovic, Doob and Sachs \cite{cds} and Feng and Li \cite{feng}.

Suppose that $\mathcal{X}=(V,E)$ is a simple graph (without loops and multiple edges), for a fixed vertex $v$ in $V$, let $\Gamma(v)$ denote the set of
neighbors of $v$, that is, the set of all vertices adjacent to $v$. The number of neighbors of $v$ is called the valency of $v$. For a subset
$\Omega \subseteq V $, its neighborhood is defined as the set
$$\Gamma(\Omega) = \{u \in V \setminus \Omega : u \mbox{ is adjacent to some vertex } v \in  \Omega\}.$$
 If the valencies of all vertices are equal, say $k$, then  $\mathcal{X}$ is called an $k$-{\it regular graph}.
 {\defn An $(|V|, k, c)$-expander is a $k$-regular graph $\mathcal{X }= (V,E)$ on $|V|$ vertices such
that every subset $\Omega \subseteq V$ with $|\Omega| \leq |V|/2$ is connected by edges of $\mathcal{X }$ to at
least $c|\Omega|$ vertices outside the set $\Omega$, that is,
\begin{equation*}
   |\Gamma(\Omega)| \geq c|\Omega| \mbox{ for all } \Omega \subseteq V \mbox{ with }|\Omega|\leq |V|/2.
\end{equation*}}
\ \  Interestingly, the difference $k-\lambda_2$ (known as the spectral gap) between
the degree $k$ of a $k$-regular graph and its second-largest eigenvalue $\lambda_2$
gives us a lot of information about the expansion properties of the corresponding
graphs. For  example, the larger this difference is, the better expansion properties the
graph has since we have the following {\it Expander Crossing Lemma.}
{\lem (\cite{jukna} Proposition 15.2) Let $\mathcal{X}=(V,E)$ be an $k$-regular graph with the second largest eigenvalue $\lambda_2$ and let $\mathcal{X}=\Omega_1\cup \Omega_2$ be a partition of $\mathcal{X}$. Then
$$e(\Omega_1,\Omega_2)\geq \frac{(k-\lambda_2)|\Omega_1|\cdot|\Omega_2|}{|V|},$$
where $e(\Omega_1,\Omega_2)$ is the number of edges between $\Omega_1$ and $\Omega_2$.}

From this lemma, we know that the bigger the spectral gap in which a communication network has, the less bottlenecks hide in it.

{\defn A connected $k$-regular graph is called a Ramanujan graph if $\lambda\leq 2\sqrt{k-1}$ for all
eigenvalues $\lambda$ of $\mathcal{X}$ with $|\lambda|\neq k$.}

Examples of Ramanujan graphs include the clique, the biclique, and the Petersen graph. As Murty's survey paper \cite{murty} notes, Ramanujan graphs ``fuse diverse branches of pure mathematics, namely, number theory, representation theory, and algebraic geometry". As an example of this, a regular graph is Ramanujan if and only if its Ihara zeta function satisfies an analog of the Riemann hypothesis \cite{murty}. Ramanujan graphs have good expansion properties. However, the known constructions of infinite families of Ramanujan graphs are rare. The only known integers $k$ for the existence of infinite families of $k$-regular Ramanujan graphs are:

(1) $k=p+1$, $p$ an odd prime, \cite{lps};

(2) $k=3$, \cite{chiu};

(3) $k=q+1$, $q$ a prime power, \cite{mor}.

So an open challenging problem is to find more infinite families of Ramanujan graphs. Particularly, ``can we construct $k$-regular Ramanujan graphs for every positive integer $k$? Currently no
one seems to know. $k = 7$ is the first unknown case," see {\rm http://www.cs.huji.ac.il/$^\sim$nati/PAPERS/nips04$\underset{-}{}$talk.pdf}.

{\defn A $k$-regular graph is called a $(|V|, k,\lambda,\mu)$-{\it strongly regular graph} (in short srg) if for every pair of vertices $u,v\in V$, the number of the common neighbors of $u$ and $v$ is
\begin{equation*}
    \left\{\begin{array}{ll}
             \lambda, & \mbox{ if $u$ and $v$ are adjacent}, \\
             \mu, & \mbox{ if  $u$ and $v$ are nonadjacent}.
           \end{array}
    \right.
\end{equation*}}

Strongly regular graphs are also closely related to two-weight
codes, two-intersection sets in finite geometry, and partial difference sets. For these connections, we refer the reader to \cite{cal, ma}.
It is known that one can construct  strongly regular graphs by using partial difference sets and cyclotomic cosets, see for example, \cite{fengtao,ma}. In \cite{cao},
we defined the so called {\it generalized difference sets} (in short, GDS) as
follows.

{\defn \cite{cao} Let $C$ be a $k(>1)$-subset of an Abelian group $G$, $S$ be a
subset of $G$, $S\neq \emptyset$. Suppose that for every element
$g(\neq 1_G)$,
\begin{eqnarray*}\mu_g=\left\{
\begin{array}{ll}
  \mu_1, & \mbox{ if \ }g\in S \\
  \mu_2, & \mbox{ if \ }g\notin S,\\
\end{array} \right. \end{eqnarray*}
where $\mu_g$ is defined by the number of unordered pairs $(c_1,c_2)\in C\times C$ such that $g=c_1c_2^{(-1)}$.
  Then
$C$ is called a $(n,|S|,k,\mu_1,\mu_2)$ generalized
difference set related to $S$. If $G$ is a cyclic group, then $C$ is called a cyclic generalized
difference set (in short. CGDS).}
For an integer $t$, let $C(t)$ be the image of $C$ under the map $c\mapsto c^t$. If $C(t)$ is a translation of $C$, then $t$ is called a multiplier of $C$.

GDS is a generalization of difference set. For example.

(a) If $S$ is a subgroup of $G$, then $C$ is a group divisible
difference set of $G$ \cite{pott};

(b) If $S$ is a subgroup of $G$, $\mu_1=0$, then $C$ is a
relative difference set of $G$ \cite{pott}. Furthermore, if
$S=\{1_G\}$, then $C$ is a difference set of $G$ and denoted by $(n,k,\mu)$ the parameters of $C$;

(c) If $S=C$, then $C$ is a partial difference set of $G$ \cite{ma};

(d) If $S$ is a subgroup of $G$, and $|\mu_1-\mu_2|=1$,
then $C$ is called an almost difference set, see for example,
\cite{davis, arasu}.

Difference sets and almost difference sets have many applications in
coding theory and cryptography, the reader is referred to
\cite{jungnickel, pott,arasu} for more details.

Let $D$ be a subset of $G$. Then $D$ corresponds uniquely to an
element of the group-ring $\mathbf{Z}[G]$, namely $D\leftrightarrow
\sum_{d\in D}d$. We denote $\sum_{d\in D}d$ simply by $D$, and
$\sum_{d\in D}d^{-1}$ by $D^{(-1)}$. It is easy to prove the
following theorem.

{\thm \cite{cao}\label{thm-1} An $k$-subset $C$ of an Abelian group $G$ is a $(v,|S|,k,\mu_1,\mu_2)$-GDS related to $S$ if and only if the following
conditions hold in $\mathbf{Z}[G]$.

(1) If $S$ contains the identity of $G$, then
\begin{eqnarray}CC^{(-1)}=(k-\mu_1)1_G+\mu_1S+\mu_2(G-S);\label{7f-1}\end{eqnarray}

(2) If $S$ does not contain the identity of $G$, then
\begin{eqnarray}CC^{(-1)}=(k-\mu_2)1_G+\mu_1S+\mu_2(G-S).\label{7f-1'}  \end{eqnarray}}



In this paper, we aim to provide more expander graphs and strongly regular graphs. We will use difference sets and generalized difference sets to construct more expander graphs, strongly regular graphs and some infinite families of Ramanujan graphs in Section 2. In Section 3 and Section 4, we will construct three infinite families of Ramanujan graphs. An example is presented to illustrate the efficiency of our method. The example shows that one can obtain a strongly regular graph which is also a Ramanujan graph, more over, the constructed graph may be distance transitive.
It is also shown that for every integer $t>2$, there exists an $k$-regular Ramanujan graph with $k=2t^2-2$. It is also proved that for every odd number $m$, there is an $(2^{2m-2}+2^{m-1})$-regular Ramanujan graph. These results partially solved the long hanging open question for asking the existence of $k$-regular Ramanujan graphs for every positive integer $k$. An algorithm is presented for finding Cayley graphs which are also Ramanujan graphs by searing suitable subset in the groups.

\section{strongly regular graphs and expander graphs arising from GDSs}
Let $\mathcal{X}=(V,E)$ be a simple graph and let $A$ be its adjacency matrix. Let $\Lambda(\mathcal{X})$ be the set of eigenvalues of $A$.
For the regular graph, one has the following well known facts:
{\lem \cite{jukna} Suppose that $\mathcal{X}$ is $k$-regular.  Then

(1) $k$ is a eigenvalue of $A$ with the associated eigenvector ${\bf 1}_n=(1,1,\cdots,1)^t$;

(2) The number of components of $\mathcal{X}=(V,E)$ is equal to the algebraic multiplicity of the eigenvalue $k$ of the characteristic polynomial of $\mathcal{X}=(V,E)$;

(3) If $\mathcal{X}=(V,E)$ is connected,
then $\mathcal{X}$ is a complete graph if and only if $|\Lambda(\mathcal{X})|=2$;

(4) If $\mathcal{X}=(V,E)$ is connected,
then $\mathcal{X}$ is a srg if and only if $|\Lambda(\mathcal{X})|\leq 3$;

(5) All eigenvalues of $\mathcal{X}=(V,E)$ are real and $\lambda_1=k\geq \lambda_2\geq \cdots \geq \lambda_n\geq -k$;

(6) $\mathcal{X}=(V,E)$ is a bigraph if and only if the eigenvalues of $\mathcal{X}=(V,E)$ are symmetric about zero;

(7) $\lambda_s=\max_{||x||=1,x\perp U, \dim(U)=s-1}x^tAx$ for all $s,1\leq s\leq n$.

(8) $\mathcal{X}$ is connected if and only if $\lambda_2<k$.}

Notice that graphs with few distinct eigenvalues have some interesting applications in communication networks and theoretic computer sciences, see for example \cite{cds}.

In what follows, we use GDSs to construct some graphs with some good expanding properties.

Let $C$ be a $(n,|S|,k,\mu_1,\mu_2)$-GDS in an Abelian group $G$ of order $n$ related to $S$ and suppose that $C^{(-1)}=C$. Then we define a Cayley graph $\mathfrak{C}(G;C,S)$ as the graph with vertices set $G$ and two points $u$ and $v$ are adjacent if and only if $uv^{-1}\in C$. Define $\mathbb{C}^G$ be the set of all complex functions from $G$ to the complex field $\mathbb{C}$. For $f,g\in \mathbb{C}^G$, define the inner product of $f$ and $g$ by
\begin{equation*}
    (f,g)=\frac{1}{n}\sum_{x\in G}f(x)\overline{g(x)}.
\end{equation*}
Then $\mathbb{C}^G$ is a unitary space with the set of characters of $G$ as a standard orthogonal basis. For any complex matrix $M=(m_{ij})$ of size $n\times n$, one can view $f$ as a $|G|$-dimensional vector with its $i$-th coordinate $f(v_i)$, and view $M$ as an operator on $\mathbb{C}^G$ in the usual way. To say it explicitly, for any function $f$ in $\mathbb{C}^G$, we define
  \begin{equation*}
   M(f)=(\sum_{t=1}^{n}m_{1t}f(v_t),\cdots, \sum_{t=1}^{n}m_{nt}f(v_t))^T.
\end{equation*}
Thus, we have $M(f)(v_i)=\sum_{j=1}^{n}m_{ij}f(v_j)$. If $M$ is the adjacent matrix of $\mathfrak{C}(G;C,S)$, then
\begin{equation*}
    M(f)(x)=\sum_{y\in G:x\rightarrow y}f(y).
\end{equation*}
Thus for every characters $\chi$ of $G$, we have
\begin{eqnarray*}
  &&M(\chi)(x) = \sum_{y:x\rightarrow y}\chi(y)=\sum_{c\in C}\chi(xc)=\chi(x)\sum_{c\in C}\chi(c)=\chi(x)\chi(C),
\end{eqnarray*}
here we define $\chi(C)=\sum_{c\in C}\chi(c)$.
Thus, the set of eigenvalues of $\mathfrak{C}(G;C,S)$ is $\{\chi(C):\chi\in \widehat{G}\}$, where $\widehat{G}$ is the character group of $G$. Since $M$ is symmetric, $\chi(C)$ is a real number for all $\chi\in \widehat{G}$.
Moreover, by Theorem \ref{thm-1}, we know that $\chi(S)$ is a a real number for all $\chi\in \widehat{G}$, and
\begin{equation*}
    \chi(C)=\left\{\begin{array}{ll}
                    \pm \sqrt{k-\mu_1+(\mu_1-\mu_2)\chi(S)}, & \mbox{ if $1_G\in S$}, \\
                    \pm  \sqrt{k-\mu_2+(\mu_1-\mu_2)\chi(S)}, & \mbox{ if $1_G\not\in S$.}
                   \end{array}
    \right.
\end{equation*}

In the next sequel, we always assume that $1_G\in S$, and the case of $1_G\not\in S$ can be discussed similarly.

By \cite{jukna}, Proposition 15.3, also Lemma 2.1 (8), we know that a $k$-regular graph with second-largest eigenvalue $\lambda_2$ is
connected if and only if $\lambda_2 < k$.
It is obvious that the second largest eigenvalue of $\mathfrak{C}(G;C,S)$ is
\begin{equation*}
    \lambda_2=\max_{\mbox{ non-trivial }\chi\in \widehat{G} }\sqrt{k-\mu_1+(\mu_1-\mu_2)\chi(S)}
\end{equation*}
Thus we know that $\mathfrak{C}(G;C,S)$ is connected if  $-\mu_1+(\mu_1-\mu_2)|S|<k^2-k$. If $-\mu_1+(\mu_1-\mu_2)|S|<3k-4$, then $\mathfrak{C}(G;C,S)$ is a Ramanujan graph. In summary, we have
{\thm Let $C$ be a $(n,|S|,k,\mu_1,\mu_2)$-GDS in an Abelian group $G$ related to a subset $S$ with $C^{(-1)}=C$. Suppose that $\mathfrak{C}(G;C,S)$
is the corresponding Cayley graph, then all nontrivial eigenvalues of $\mathfrak{C}(G;C,S)$ are in the set $$\{ \pm\sqrt{k-\mu_1+(\mu_1-\mu_2)\chi(S)},\chi\in \widehat{G} \mbox{ is non-principal}\}.$$
If $C=S$, then $\mathfrak{C}(G;C,S)$ is a strongly regular graph.
Moreover,  $\mathfrak{C}(G;C,S)$ is connected if $-\mu_1+(\mu_1-\mu_2)|S|<k^2-k$. If $-\mu_1+(\mu_1-\mu_2)|S|<3k-4$, then $\mathfrak{C}(G;C,S)$ is a Ramanujan graph.}

When $S$ is a subgroup of $G$, we know that for every nontrivial character $\chi$ of $G$, it holds that
\begin{equation*}
    \chi(S)=\left\{\begin{array}{ll}
                     |S|, & \mbox{ if $\chi$ is principal on $S$,} \\
                     0, & \mbox{ otherwise.}
                   \end{array}
    \right.
\end{equation*}
Thus, in this case, we have
\begin{equation}\label{f-2}
    \chi(C)=\left\{\begin{array}{ll}\pm \sqrt{k-\mu_1+(\mu_1-\mu_2)|S|},& \mbox{ if $\chi$ is principal on $S$,} \\
    \pm \sqrt{k-\mu_1}, & \mbox{ otherwise.}
                   \end{array}
    \right.
\end{equation}
Thus the second largest eigenvalue of $\mathfrak{C}(G;C,S)$ is
\begin{equation*}
    \lambda_2\leq\left\{\begin{array}{ll}
                       \sqrt{k-\mu_1+(\mu_1-\mu_2)|S|}, & \mbox{ if $\mu_1>\mu_2$,} \\
                       \sqrt{k-\mu_1}, & \mbox{ if $\mu_1<\mu_2$.}
                     \end{array}
    \right.
\end{equation*}
 Hence, we have the following result.
{\thm Let $C$ be a $(n,|S|,k,\mu_1,\mu_2)$-GDS in an Abelian group $G$ related to a subgroup $S$ with $C^{(-1)}=C$. Suppose that $\mathfrak{C}(G;C,S)$
is the corresponding Cayley graph, then all nontrivial eigenvalues of $\mathfrak{C}(G;C,S)$ are in the set $$\{ \pm\sqrt{k-\mu_1+(\mu_1-\mu_2)|S|}, \pm\sqrt{k-\mu_1}\}.$$
Moreover, $\mathfrak{C}(G;C,S)$ is connected if and only if $(\mu_1-\mu_2)|S|-\mu_1<k^2-k$. Furthermore, if $(\mu_1-\mu_2)|S|\leq 3k+\mu_1-4$, then $\mathfrak{C}(G;C,S)$ is a Ramanujan graph. Particularly, if $\mu_1<\mu_2, k>1$, then $\mathfrak{C}(G;C,S)$ is always a Ramanujan graph.}

It should be noted that there are plenty of GDSs in some Abelian groups in general. Below, we provide an algorithm for finding the GDSs in the cyclic group $\mathbb{Z}_n=\mathbb{Z}/(n)$.

\noindent{\bf Algorithm 1}:  An algorithm for finding the GDSs in the cyclic group
$\mathbb{Z}_n$

Input: a positive integer $n$; Output: $c(x)$ and $S$.

Step 1. For every positive integer $s$, $1<s<2^n-1$, express $s$ as
the $2-$adic number, $s=\sum_{i=0}^{n-1}a_i2^i$;

Step 2. Generate the Hall polynomial
$$c(x)=\sum_{i=0}^{n-1}a_ix^i;$$

Step 3. Calculate $c(x)c(x^{n-1})$ mod $(x^n-1)$;

Step 4. If there are two integers $\mu_1,\mu_2$ such that
$$c(x)c(x^{n-1})\equiv(k-\mu_1)1_G+(\mu_1-\mu_2)S+\mu_2G \quad \mod (x^n-1),$$
where $k=\sum_{i=0}^{n-1}a_i, S$ corresponds to a subset of $G$, then print
$c(x)$ and $S$;

Step 5. If $s<2^n-1$, then go to Step 1;

Step 6. END.

Now we provide some examples to illustrate that our method works well.

\noindent {\bf Example 1} Let $n=20$ and $G=\mathbb{Z}_{n}$. Let $C=\{4,8,12,16\}$ and $c(x)=x^4+x^8+x^{12}+x^{16}$. Then it is easy to verify that
\begin{equation*}
    c(x)^2\equiv 4+3(x^4+x^8+x^{12}+x^{16})\ ({\rm mod}\ x^n-1).
\end{equation*}
Thus $C$ is a $(20,16,4,0,3)$-GDS in $G$, and the corresponding Cayley graph $\mathfrak{C}(G;C,S)$ has $4$ components, each component is a complete graph. The spectral of the graph are $4$ with multiplicity $4$ and $-1$ with multiplicity $16$.
If we put $C'=\{2,6,14,18\}$, then $C'$ is also a $(20,16,0,3)$-GDS in $G$. In this case, the corresponding Cayley graph $\mathfrak{C}(G;C',S')$ is a bipartite graph, The spectral of the graph are $\pm 4$ with multiplicity $2$ and $\pm 1$ with multiplicity $8$. Interestingly, if we put $C=\{3,4,8,12,16,17\}$, then the corresponding Cayley graph $\mathfrak{C}(G;C,S)$ is a Ramanujan graph.

Below, we also provide an algorithm for finding a subset $C$ in the cyclic group $\mathbb{Z}_n$ such that the Cayley graph $\mathfrak(\mathbb{Z}_n,C)$ is a Ramanujan graph.

\noindent{\bf Algorithm 2:} An algorithm for finding Ramanujan graphs from Cayley graph corresponding to a subset $C$ in the cyclic group $\mathbb{Z}_n$.

Input: a positive integer $n$; Output: the subset $C$.

Step 1. For every positive integer $s$, $1<s<2^{\lceil n/2-1\rceil}$, express $s$ as
the $2-$adic number, $s=\sum_{i=0}^{\lceil n/2-1\rceil}a_i2^i$, denote $wt(s)=\sum_{i=0}^{\lceil n/2-1\rceil}a_i$;

Step 2. Generate the Hall polynomial
$$c(x)=\sum_{i=0}^{\lceil n/2-1\rceil}a_i(x^i+x^{n-i});$$

Step 3. For every integer $a$, $1\leq a\leq n-1$, calculate $C_a:=c(\xi_n^a)$, where $\xi_n$ is a primitive $n$-th root of unity;

Step 4. If $C_a> 2\sqrt{wt(s)-1}$, then go to Step 1;

Step 5. print $C=\{i,n-i: a_i\neq 0, 0\leq i\leq \lceil n/2-1\rceil\}$;

Step 6. END.

Using this algorithm, we obtain the following Table 1 of pairs $(n,s)$ which correspond to the subset in the cyclic group $\mathbb{Z}_n, 15\leq n\leq 30$, the associated Cayley graph related to each subset is a Ramanujan graph.

\begin{center}
Table 1: Some pairs $(n,s)$ which lead to Ramanujan graphs
\begin{tabular}{|c|c||c|c|}
  \hline
$ n $& $s$ & $n$ &$ s$ \\
 \hline
  15& 175,111,191& 20& 973,1007\\
    \hline
     21 & 1535,959,1471 & 27 &12031 ,12287 , 7935 \\
  \hline
   24 & 2030,2045,3935,1982,4055,959,1007&28&3055, 7930, 15773,16319,8127\\
  &4015,4063,3967,4087,3835,2039,3959&&\\
  &3581,1919,3070,1535,2043,1791,2015&&\\
  \hline
  30&32491,15871,16367,31979,32751&&\\
   &32255,32239,32507 &&\\
  \hline
\end{tabular}

\end{center}
{\small Note that for each number $s$ in the table, we expand it in $2$-adic numbers, for example, for $n=20$, and $s=973$, we have $973=1011001111_2$, thus the corresponding subset is $C=\{1,3,4,7,8,9,11,12,13,16,17,19\}$ (remove the number $10$ since $\xi^{10}=-1$ for any primitive $20$th root of unity $\xi$). It can be verified that $C$ is a reversible difference set in $G$, and $\mathfrak{C}(G,C)$ is a Ramanujan graph.}

It is obvious that if $C$ is actually a difference set in $G$, then as a consequence of Theorem 2.2 and Theorem 2.3, we have the following result.

{\cor If $C$ is a difference set in an Abelian group $G$ with minus one as a multiplier, then the associated Cayley graph of $C$ (or one of the translations of $C$) is a Ramanujan graph.}

Recall that there are some difference sets both in Abelian groups and non-Abelian groups with minus one as a multiplier. For example, we have the following result due to Kraemer \cite{kra}, Jedwab \cite{jedwab} and Turyn \cite{turyn}:
{\lem \label{lem-2.4}An Abelian $2$-group $G$ of order $2^{2 t+2}$ has a Hadamard
difference set, i.e., with parameters $(2^{2 t+2},2^{2 t+1}\pm 2^t, 2^{2 t}\pm 2^t)$, if and only if the exponent of the group is less than or equal to $2^{t+2}$.}

It is obvious that if $C$ is a difference set in an elementary abelian $2$-group, then $-1$ is a multiplier of $C$. Thus by Lemma 2.5 and Theorem 2.2, we have the following infinite family of Ramanujan graphs.

{\prop For every positive integer $t$, there is a Hadamard
difference set in the group $\mathbb{Z}_2^t$ (the direct product of $t$ copies of $\mathbb{Z}_2$), and the
 corresponding Cayley graph is a Ramanujan graph.}

Let $q$ be a prime power and let $s$ be a positive integer. In \cite{mcfarland}, McFarland  constructed
difference sets with the parameters
\begin{equation}\label{f-524}
    \left(q^{s+1}\left(\frac{q^{s+1}-1}{q-1}+1\right), q^s\left(\frac{q^{s+1}-1}{q-1}\right),q^s\left(\frac{q^{s}-1}{q-1}\right)\right)
\end{equation}
in a group $G$ which is the direct product
$G=E\times K,$
where $E$ is an elementary Abelian group of order $q^{s+l}$ and $K$ is any group
of order $\frac{q^{s+1}-1}{q-1}+1$. He also provided a sufficient and necessary conditions for $-1$ being a multiplier of this difference set.
Thus by Lemma 2.4 and Theorem 2.1 again, we have

{\prop Let $q$ be a prime power and let $s$ be any positive integer. If $E$ is an elementary Abelian group of order $q^{s+l}$ and $K$ is a suitable group
of order $\frac{q^{s+1}-1}{q-1}+1$, then there is a difference set $C$ in the group $E\times K$ with $-1$ as a multiplier. As a consequence, the Cayley graph associated to $C$ (or a translation of $C$) is a Ramanujan graph.}

For any even number $s\geq 2$, let $\mathbb{F}_{3^{s+1}}$ be the finite field with $3^{s+1}$ elements and let $\mathbb{F}_{3^{s+1}}^*$ be the multiplicative group of $\mathbb{F}_{3^{s+1}}$. Then since $s$ is even, $\mathbb{F}_{3^{s+1}}^*=\mathbb{Z}_2\times K$ for a group $K$ of order $\frac{3^{s+1}-1}{2}$. Let $G=E\times K$ where $E$ is additive group of $\mathbb{F}_{3^{s+1}}$. In 1979, M. Miyamoto proved that there exists a difference set in the group $G$ with parameter
\begin{equation}\label{f-5242}
    \left(\frac{3^{s+1}(3^{s+1}-1)}{2}, \frac{3^s(3^s+1)}{2}, \frac{3^s+1}{2}\right)
\end{equation}
which has $-1$ as a multiplier, see http://hdl.handle.net/2433/104316. Thus, we know that

{\prop For any even positive integer $s$, there is a $k$-regular Ramanujan graph, where $k=\frac{3^s(3^s+1)}{2}$.}

Note that in \cite{wilson}, an infinite family of reversible Hadamard difference sets were presented, thus, one can also construct some Ramanujan graphs using such Hadamard difference sets. We omit the details.
\section{An elementary construction of an infinite family of Ramanujan graphs}
In this section, we will provide a construction of an infinite family of Ramanujan graphs.



Let $E=\mathbb{Z}_{s}$ and $K=\mathbb{Z}_r$ be the cyclic group of order $s$ and order $r$ respctively, where $s,r\geq 4$ is an even number. Let $G=E\times K$ and let $C=\{x_1,\cdots,x_{s'}\}$ which forms a proper subgroup of $E$ and let
$\mathfrak{C}(G,C)$ be the Cayley graph corresponding to the set $C$. In order to remove the loops, we should delete the element zero from $C$. Let $C_0=C\setminus \{0\}$. Then $\mathfrak{C}(G,C_0)$ has $sr/{s'}$ connected components, each one is a complete graph. Let $A$ be the adjacent matrix of $\mathfrak{C}(G,C_0)$, then there exists a permutation matrix $P$ such that
\begin{equation}\label{f-5245}
   P^{-1}AP={\rm diag}[J_1,J_1,\cdots,J_1 ]
\end{equation}
where $J_1=J-I$.
Hence the eigenvalues of $\mathfrak{C}(G,C_0)$ are $s'-1$ with multiplicity $sr/{s'}$ and $-1$ with multiplicity $sr-sr/s'$.
Indeed, for each non-trivial character $\chi$ of $G$, we have
\begin{equation}\label{f-5244}
    \chi(C_0)=\left\{\begin{array}{ll}
                       s'-1,  & \mbox{ if $\chi$ is trivial on $C$,} \\
                       -1, & \mbox{otherwise}.
                     \end{array}
    \right.
\end{equation}

Putting $C_1=\{2a|a\in \mathbb{Z}_r\}\setminus\{0\}$, and $$D=[(C_0\times K)\cup (E\times C_1)]\setminus [(C_0\times K)\cap (E\times C_1)] \subset G,$$ we can show that
the Cayley graph $\mathfrak{C}(G,D)$ is a Ramanujan graph under the some hypothesis. That is:

{\thm \label{thm-33} The above constructed graph $\mathcal{X}=\mathfrak{C}(G,D)$ is a Ramanujan graph if $s\geq 4$ and $2s>r\geq s$ or if $r\geq 4$ and $2r>s\geq r$.}

\begin{proof}Without lose of generality, we assume that $r\geq s$. It is obvious that $D=A\cup B$ where $A=(C_0\times K)\setminus (C_0\times C_1)$, $B=(E\times C_1)\setminus (C_0\times C_1)$. Since $A\cap B=\emptyset$, in order to compute the eigenvalues of $\mathfrak{C}(G,D)$, we only need to evaluated the exponential sums on $A,B$. For any $\chi\in \widehat{G}$, there exists $\chi_1\in \widehat{E}$ and $\chi_2 \in \widehat{K}$ such that $\chi=\chi_1\chi_2$. Thus we have
\begin{equation*}
    \chi(D)=\sum_{c\in D}\chi(c)=\chi_1(C_0)(\chi_2(K)-\chi_2(C_1))+(\chi_1(E)-\chi_1(C_0))\chi_2(C_1).
\end{equation*}
If $\chi$ is the principal character of $G$, then $\chi(D)=sr/2-2$. Let $\chi$ be a nontrivial character of $G$.

 (1) If $\chi_1$ is the principal, $\chi_2$ is not the principal and trivial on $C_1$, then
 \begin{equation*}
   \chi(D)=(s'-1)(0-\frac{r-2}{2})+(s-(s'-1))\frac{r-2}{2}=r-2.
 \end{equation*}

 (2) If $\chi_1$ is the principal, $\chi_2$ is nontrivial on $C_1$, then
 \begin{equation*}
   \chi(D)=(s'-1)(0-(-1))+(s-(s'-1))(-1)=-s+2s'-2=-2.
 \end{equation*}

(3) If $\chi_1$ is not the principal and trivial on $C_0$, and $\chi_2$ the principal, then
\begin{equation*}
    \chi(D)=(s'-1)(r-\frac{r-2}{2})+(0-(s'-1))\frac{r-2}{2} =s-2.
\end{equation*}

(4) If $\chi_1$ is nontrivial on $C_0$, and $\chi_2$ the principal, then
 \begin{equation*}
   \chi(D)=(-1)\frac{r+2}{2}+(0-(-1))\frac{r-2}{2}=-2.
 \end{equation*}

 (5) If $\chi_1$ is not the principal and trivial on $C_0$, and $\chi_2$ the not principal and trivial on $C_1$, then
 \begin{equation*}
   \chi(D)=(s'-1)(0-\frac{r-3}{2})+(0-(s'-1))\frac{r-3}{2}=-(s-2)(r-2)/4.
 \end{equation*}

 (6) If $\chi_1$ is nontrivial on $C_0$, and $\chi_2$  is not the principal and trivial on $C_1$, then
\begin{equation*}
     \chi(D)=(-1)(0-\frac{r-2}{2})+(0-(-1))\frac{r-2}{2})=r-2.
\end{equation*}

 (7) If $\chi_1$ is not the principal but trivial on $C_0$, and $\chi_2$  is nontrivial on $C_1$, then
 \begin{equation*}
    \chi(D)=(s'-1)(0-(-1))+(0-(s'-1))(-1)=2s'-2=s-2.
 \end{equation*}

 (8) If $\chi_1$ is nontrivial on $C_0$, and $\chi_2$ the nontrivial on $C_1$, then
 \begin{equation*}
    \chi(D)=-1(0-(-1))+(0-(-1))(-1)=-2.
 \end{equation*}
 Thus, we know that $\mathcal{X}$ is a Ramanujan graph if and only if $\max\{r-2,s-2\}<2\sqrt{sr/2-3}$.
Now, it is easy to see that
\begin{equation}\label{f-271}
    r-2<2\sqrt{sr/2-3} \mbox{ if and only if } (2s+4-r)r>16,
\end{equation}
and
\begin{equation}\label{f-272}
     s-2<2\sqrt{sr/2-3} \mbox{ if and only if } (2r+4-s)s>16.
\end{equation}
Thus, if $s\geq 4$ and $2s>r\geq s$, then (\ref{f-271}) and (\ref{f-272}) are valid simultaneously. Therefore,
The desired result follows with the hypothesis. \end{proof}

{\rem (1) Our construction is similar to Menon's method (see \cite{menon} ) to produce Hadamard difference sets. However, the set $D$ in above theorem is not a Hadamard difference set in general.

(2) Theorem \ref{thm-33} indicates that for every even integer $k>4$, we can find even numbers $r$ and $s$ satisfying $rs=2(k+2)$, then the resulting graph is an $k$-regular Ramanujan graph if $s\geq 4$ and $2s>r\geq s$ or if $r\geq 4$ and $2r>s\geq r$. Particularly, Letting $r=s=2t^2$, we get that there exists an $k$-regular Ramanujan graph with $k=2t^2-2$ for every integer $t$. }

Below, we give an example to illustrate the efficiency of our method.

\noindent{\bf Example 2} Let $s=r=4$. Then $E=K=\mathbb{Z}_4$. In this case, we have
$C_0=\{2\}$, $C_1=\{2\}$, and
\begin{equation*}
    D=\{(2,0),(0,2),(1,2),(2,1),(3,2),(2,3)\}.
\end{equation*}
In order to present the graph, we label the elements of $G=E\times K$ as follows:
\begin{center}

\mbox{Table 2. Labels of elements of $G$}.

      \begin{tabular}{|ccccccccc|}
     \hline
     \mbox{label}&1 & 2 & 3 & 4 &5 & 6 & 7 & 8  \\
     \hline
     \mbox{elements}&(0,1) & (0,2)& (0,3) & (1,0) & (1,1) & (1,2) & (1,3) & (2,0)\\
     \hline
     \mbox{label}& 9 &10 & 11 &12 & 13 & 14 & 15 & 16 \\
      \hline
     \mbox{elements}& (2,1) & (2,2) & (2,3)& (3,0) & (3,1) & (3,2) & (3,3) & (0,0) \\
     \hline
   \end{tabular}
    \end{center}
The associated Cayley graph $\mathfrak{C}(G,D)$ has the vertices as the set $G$ and the neighbors of each point are listed in the following Table 3:
\begin{center}

\mbox{Table 3. Neighbors of the vertices}.

   \begin{tabular}{|c|c||c|c|}
     \hline
     \mbox{Vertex}&\mbox{neighbors}& \mbox{Vertex}&\mbox{neighbors} \\
     \hline
     \mbox{1}&3, 7, 8, 9, 10, 15 & \mbox{9}&1, 2, 7, 11, 15, 16 \\
     \hline
     \mbox{2}&4, 9, 10, 11, 12, 16& \mbox{10}&1, 2, 3, 4, 8, 12  \\
      \hline
     \mbox{3}&1, 5, 8, 10, 11, 13& \mbox{11}&2, 3, 5, 9, 13, 16 \\
     \hline
      \mbox{4}&2, 6, 10, 12, 13, 15& \mbox{12}& 2, 4, 5,7, 10, 14  \\
      \hline
      \mbox{5}&3, 7, 11, 12, 13, 14& \mbox{13}& 3, 4, 5, 6, 11, 15 \\
      \hline
      \mbox{6}&4, 8, 13, 14, 15, 16& \mbox{14}&5, 6, 7, 8, 12, 16 \\
      \hline
      \mbox{7}&1, 5, 9, 12, 14, 15 & \mbox{15}&1, 4, 6, 7, 9, 13 \\
      \hline
      \mbox{8}&1, 3, 6, 10, 14, 16 & \mbox{16}&2, 6, 8, 9, 11, 14 \\
      \hline
   \end{tabular}
    \end{center}
$\mathfrak{C}(G,D)$ has the spectral $6$ with multiplicity $1$ and $2$ with multiplicity $6$ and $-2$ multiplicity $9$. This is a srg and a Ramanujan graph! The automorphism group of $G$ has order $1152=7^2\cdot 3^2$ and has the following set of generators:
\begin{eqnarray*}
     (4, 12)(5, 13)(6, 14)(7, 15),\\
    (2, 4)(6, 16)(9, 15)(11, 13),\\
    (2, 6)(4, 16)(8, 10)(9, 15)(11, 13)(12, 14),\\
    (3, 7)(4, 16)(8, 15)(9, 10)(11, 12)(13, 14),\\
    (1, 2)(3, 4, 8, 12)(5, 13, 6, 14)(7, 11, 15, 16).
\end{eqnarray*}
The graph $\mathfrak{C}(G,D)$ is primitive, i.e., the automorphism group of $\mathfrak{C}(G,D)$ is primitive. $\mathfrak{C}(G,D)$ is also distance transitive which means that if two pairs of vertices in $\mathfrak{C}(G,D)$ share the common distance, then there exists an automorphism sending one pair to another. The  intersection array of $\mathfrak{C}(G,D)$ is $[6,3,1,2]$.
 Therefore, the graph has a very strongly symmetry. And note that the set $D$ is not constructed by the product of two Paley difference sets.
\section{two more classes of expander graphs arising from finite fields}

In this section, we will present more construction of expander graphs and investigate their expander properties.

\subsection{The first class}

Let $q=2^m$ be a power of $2$ and let $\mathbb{F}_q$ be the finite field with $q$ elements. Let ${\rm Tr}: \mathbb{F}_q\rightarrow \mathbb{F}_2$ be the trace map. Before we going to define the graphs, we
recall the definition and properties of an exponential sum-the Kloosterman sum:
\begin{eqnarray}
  \label{equ-k}k_m(a,b) &=& \sum_{x\in \mathbf{F}_{2^m}^*}(-1)^{{\rm Tr}_m(ax+bx^{-1})}, a,b\in \mathbb{F}_q.
 \end{eqnarray}
 It is easy to check that $k_m(a,b)=k_m(ab,1)=k_m(1,ab)=:k_m(ab)$. Moreover, The Kloosterman sum $k_m(a,b)$ can be calculated recursively; that is:
 if we define
 \begin{equation*}
    k_m^{(s)}(a)=\sum_{\gamma\in \mathbb{F}_{q^s}^*}\chi^{(s)}(a\gamma+\gamma^{-1}), a\in \mathbb{F}_q,
 \end{equation*}
where $\chi^{(s)}$ is the lifting of $\chi(x)=(-1)^{{\rm Tr}(x)}$ to $\mathbb{F}_{q^s}$.  Then
 \begin{equation}\label{f-41}
    k_m^{(s)}(a)=-k_{m}^{(s-1)}(a)k_m^{(1)}(a)-qk_{m}^{(s-2)}(a),
 \end{equation}
 where we put $k_m^{(0)}(a,b)=-2$ and $k_m^{(1)}(a)=k(a)$. And, for all $a,b\in \mathbb{F}_q^*$, one has that
  \begin{equation*}
   |k_m(a,b)|\leq 2\sqrt{q},
\end{equation*}
see \cite{lidl} for details.
 Using the recursive relation, one can obtain the formulae for $k_m(1)$ as follows:
 {\lem \cite{lidl} $k_m(1)=-\omega_1^m-\omega_2^m$, where $\omega_1,\omega_2$ are the (complex) roots of the equation $x^2+x+2=0$. Moreover, one has that $k_1(1)=1,k_2(1)=3$ and
  \begin{equation}\label{f-43}
    k_{m+2}(1)+k_{m+1}(1)+2k_{m}(1)=0.
\end{equation}}

 {\lem \cite{carlitz} The value of $k_m(1)$ is
\begin{eqnarray*}
  k_m(1)&=& -\sum_{j=0}^{\lfloor m/2\rfloor}(-1)^{m-j}\frac{m}{m-j}{{m-j}\choose{j}}2^j.
\end{eqnarray*}}

Note also that the values of Kloosterman sums over $\mathbb{F}_{2^m}$
were determined by Lachaud and Wolfmann in \cite{lachaud}.

\begin{lem}\label{lem-43} The set $\{k_m(\lambda),\lambda\in \mathbb{F}_{2^m}\}$ is the set of all the integers $s\equiv -1({\rm mod}\ 4)$ in the range
\begin{equation*}
    \left[-2^{\frac{m}{2}+1},2^{\frac{m}{2}+1}\right].
\end{equation*}
\end{lem}

We define a subset $D\in \mathbb{F}_q$ as
\begin{equation*}
    D=\{z\in \mathbb{F}_q| z\neq 0 \mbox{ and } {\rm Tr}(z)={\rm Tr}(z^{-1})=1\}.
\end{equation*}
Obviously, the cardinality of $D$ is
\begin{eqnarray*}
 |D| &=& \frac{1}{4}\sum_{x\in \mathbb{F}_q^*}(1-(-1)^{{\rm Tr}(x)})(1-(-1)^{{\rm Tr}(x^{-1})}) \\
  &=& \frac{1}{4}(k_m(1)-2\sum_{x\in \mathbb{F}_q^*}(-1)^{{\rm Tr}(x)}+2^m-1)\\
 &=&\frac{1}{4}(k_m(1)+2^m+1).
\end{eqnarray*}

Now, let $G$ be the additive group of $\mathbb{F}_q$. Then the Cayley graph $\mathfrak{C}(G,D)$ has $\mathbb{F}_q$ as the vertices set and two points $u,v$ are adjacent if and only if
$u-v\in D$. Following the discussion as in Section 2, we know that the eigenvalues of $\mathfrak{C}(G,D)$ are in the set
\begin{equation*}
   \{\chi(D): \chi \mbox{ is a additive character of $\mathbb{F}_q$}\}.
\end{equation*}
Every additive character of $\mathbb{F}_q$ can be written as $\chi_a(x)=(-1)^{{\rm Tr}(ax)}$ for some $a\in \mathbb{F}_q$. When $a=0$, $\chi_0$ is the trivial character. For $a\neq 0$,
\begin{eqnarray*}
 \chi_a(D)&=&\sum_{z\in D}\chi_a(z) \\
   &=&\frac{1}{4}\sum_{z\neq 0}\chi_a(z)(1-\chi(z))(1-\chi(z^{-1}))\\
   &= &\frac{1}{4}(-1-k_m(a)+k_m(a+1)-\sum_{z\neq 0}\chi_1((a+1)z))\\
   &=&\left\{\begin{array}{cc}
               -(2^m+1+k_m(1))/4,& \mbox{ if $a=1$}, \\
               (-k_m(a)+k_m(a+1))/4,& \mbox{ if $a\neq 1$}.
             \end{array}
   \right.
\end{eqnarray*}

Since $|k_m(a)|\leq 2\sqrt{q}$ for all $a\in \mathbb{F}_q$, we see that the second largest eigenvalue of the graph $\mathfrak{C}(G,D)$ is less than $\sqrt{q}$. If $k_m(1)>3$, then we find that $2\sqrt{|D|-1}>\sqrt{q}$, thus, in this case, $\mathfrak{C}(G,D)$ is a Ramanujan graph. Moreover, when $a\neq 1$, one has that $\chi_a(D)=-\chi_{a+1}(D)$, thus, the eigenvalues is symmetric about zero, and then $\mathfrak{C}(G,D)$ is a bipartite graph. In summary, we have

{\thm \label{thm-4.4} Let $G$ be the additive group of the finite field $\mathbb{F}_{2^m}$, and let $D$ be a subset in $G$ which is defined by
\begin{equation*}
    D=\{z\in \mathbb{F}_q| z\neq 0 \mbox{ and } {\rm Tr}(z)={\rm Tr}(z^{-1})=1\}.
\end{equation*}
Then the Cayley graph $\mathfrak{C}(G,D)$ is a bipartite (bigraph) graph, the valency of it is $-(2^m+1+k_m(1))/4$ and all other eigenvalues are $(-k_m(a)+k_m(a+1))/4$, $a\neq 1$. If the Kloosterman sum $k_m(1)$ is bigger than $3$, then $\mathfrak{C}(G,D)$ is a Ramanujan graph.}

We claim that the above constructed graphs form an infinite family of Ramanujan graphs. By Theorem 4.4, it is only need to
prove that there are infinite number $m$ such that $k_m(1)>3$. This is indeed the case. For otherwise, there should exist a positive integer $N$ such that $k_m(1)<3$ for all $m>N$. Then by (\ref{f-43}), we have $3>k_m(1)>-9$. By Lemma \ref{lem-43}, $k_m(1)\in \{-1,-5\}$ for all $m>N$. This is contrary to (\ref{f-43}).

Note that if we define
\begin{equation*}
    D_{i,j}=\{z\in \mathbb{F}_q| z\neq 0 \mbox{ and } {\rm Tr}(z)=i, {\rm Tr}(z^{-1})=j\}, i=0,1,j=0,1,
\end{equation*}
Then
we have
\begin{equation}\label{f-10204}
   |D_{0,0}|=2^{m-2}+(k_m-3)/4,  |D_{1,0}|=|D_{0,1}|=2^{m-2}-(k_m+1)/4.
\end{equation}
These identities come from the following computation:
\begin{eqnarray*}
  |D_{i,j}| &=& \frac{1}{4}\sum_{z\in \mathbb{F}_q^*}\left[1+(-1)^i\chi(z)\right]\left[1+(-1)^j\chi(z^{-1})\right] \\
   &=& \frac{1}{4}[2^m-1-(-1)^j-(-1)^i+(-1)^{i+j}k_m(1)].
\end{eqnarray*}
One can also use the sets $D_{i,j}$ to construct Cayley graphs. We find that $\mathfrak{C}(G,D_{1,0})$ has the same structure as $\mathfrak{C}(G,D_{1,1})$,
 but the graph $\mathfrak{C}(G,D_{0,1})$ and $\mathfrak{C}(G,D_{0,0})$ both has two connected components. All these facts can be proven by using the same method applied in the proof of Theorem \ref{thm-4.4}, we omit the details.

It is known that some bipartite graphs yield good error-correcting codes. A bipartite graph $\mathcal{X} =
(A,B,E)$ that's a good expander, yields a good error-correcting code
that comes with efficient encoding and decoding algorithms. Roughly speaking,
if $A$ has $n$ vertices, there is a one-to-one correspondence between subsets of
$A$ and $n$-bit messages. The vertices in $B$ correspond to parity check
bits. It turns out that if $\mathcal{X}$ is a good expander, then not only is the
resulting (linear) code good, it can be provably corrected by a simple
belief-propagation scheme applied to $\mathcal{X}$. Theorem 4.4 shows that the above constructed graphs are good candidates for constructing good error-correcting codes.
\subsection{The second class}
Before we going to introduce our next graph,  we need to define some notations.

Let $q=2^n$, $n=2m$ and let $\mathbb{F}_q$ be the finite field with $q$ elements. Denote the subgroup of $2^m+1$-th roots of unity in $\mathbb{F}_q$ by $\mathfrak{S}$, i.e., $\mathfrak{S}=\{z\in \mathbb{F}_q| z^{2^m+1}=1\}$. For every $x\in \mathbb{F}_q$, there is a unique polar decomposition of $x$ as $x=yz$ where $y\in \mathbb{F}_{2^m}$ and $z\in \mathfrak{S}$. Denote $x^{2^m}$ by $\overline{x}$. Then for every $x\in \mathbb{F}_q^*$, $x\in \mathbb{F}_{2^m}$ if and only if $x=\overline{x}$ and $x\in \mathfrak{S}$ if and only if $\overline{x}=x^{-1}$. it is evident that for every $x\in \mathbb{F}_q^*$, one has that $x+\overline{x}, x\overline{x}\in \mathbb{F}_{2^m}$.

We define a subset $D$ in $\mathbb{F}_q$ by
\begin{equation*}
    D=\{x\in \mathbb{F}_q^*|{\rm Tr}_{\mathbb{F}_{2^m}/\mathbb{F}_2}(x+\overline{x})={\rm Tr}_{\mathbb{F}_{2^m}/\mathbb{F}_2}(x\overline{x})=1\}.
\end{equation*}
Then the cardinality of $D$ is
\begin{eqnarray*}
  |D|&=& \frac{1}{4}\sum_{x\in \mathbb{F}_q^*}(1-(-1)^{{\rm Tr }_{\mathbb{F}_{2^m}/\mathbb{F}_2}(x+\overline{x})})(1-(-1)^{{\rm Tr}_{\mathbb{F}_{2^m}/\mathbb{F}_2}(x\overline{x})}) \\
  &=& \frac{1}{4}\left[2^n-1-U_1-U_2+W\right]
\end{eqnarray*}
where \begin{eqnarray*}
U_1 &=& \sum_{x\in \mathbb{F}_q^*}(-1)^{{\rm Tr}_{\mathbb{F}_{2^m}/\mathbb{F}_2}(x+\overline{x})},
U_2=\sum_{x\in \mathbb{F}_q^*}(-1)^{{\rm Tr}_{\mathbb{F}_{2^m}/\mathbb{F}_2}(x\overline{x})},\\
W&=&\sum_{x\in \mathbb{F}_q^*}(-1)^{{\rm Tr}_{\mathbb{F}_{2^m}/\mathbb{F}_2}(x+\overline{x}+x\overline{x})}.
\end{eqnarray*}
Let $x=yz$, $y\in \mathbb{F}_{2^m}$ and $z\in \mathfrak{S}$. Then
\begin{eqnarray*}
  U_1 &=& \sum_{y\in \mathbb{F}_{2^m}^*,z\in \mathfrak{S}}(-1)^{{\rm Tr}_{\mathbb{F}_{2^m}/\mathbb{F}_2}(y(z+z^{-1}))}.
\end{eqnarray*}
Since $z\in \mathfrak{S}$, we know that $z+z^{-1}=0$ if and only if $z=1$. Thus
\begin{equation*}
   U_1=2^m-1+(-1)2^m=-1.
\end{equation*}
For the number $U_2$, we have
\begin{equation*}
    U_2=\sum_{y\in \mathbb{F}_{2^m}^*,z\in \mathfrak{S}}(-1)^{{\rm Tr}_{\mathbb{F}_{2^m}/\mathbb{F}_2}(y^2)}=-(2^m+1).
\end{equation*}
Similarly, for the number $W$, one has that
\begin{equation*}
    W=\sum_{y\in \mathbb{F}_{2^m}^*,z\in \mathfrak{S}}(-1)^{{\rm Tr}_{\mathbb{F}_{2^m}/\mathbb{F}_2}(y(z+z^{-1}+1))}.
\end{equation*}
Since $z\in \mathfrak{S}$, $z+z^{-1}+1=0$ if and only if
\begin{equation}\label{f-4.5}
   z^2+z+1=0. \mbox{ i.e., } z^{\gcd(3,2^m+1)}=1 \mbox{ and $z\neq 1$}.
\end{equation}
When $m$ is even, (\ref{f-4.5}) has no solution. When $m$ is odd, (\ref{f-4.5}) has $2$ distinct roots. Thus, we have
\begin{equation*}
    W=\left\{\begin{array}{cc}
              (-1)(2^m+1), & \mbox{ if $m$ is even}, \\
               2(2^m-1)+(-1)(2^m-1),& \mbox{ if $m$ is odd}.
             \end{array}
    \right.
\end{equation*}
Therefore, we have
\begin{equation*}
    |D|=\left\{\begin{array}{cc}
                 2^{n-2}, & \mbox{ if $m$ is even}, \\
                 2^{n-2}+2^{m-1},& \mbox{ if $m$ is odd}.
               \end{array}
    \right.
\end{equation*}
For every nontrivial additive character $\chi_a$ of $\mathbb{F}_q$, we have
\begin{eqnarray*}
 &&\chi_a(D)=\sum_{x\in D}\chi_a(x) \\
  &=& \frac{1}{4}\sum_{x\in \mathbb{F}_q^*}\chi_a(x)(1-(-1)^{{\rm Tr }_{\mathbb{F}_{2^m}/\mathbb{F}_2}(x+\overline{x})})(1-(-1)^{{\rm Tr}_{\mathbb{F}_{2^m}/\mathbb{F}_2}(x\overline{x})})\\
  &=&\frac{1}{4}\sum_{x\in \mathbb{F}_q^*}(-1)^{{\rm Tr }_{\mathbb{F}_{2^m}/\mathbb{F}_2}(ax+\overline{ax})}(1-(-1)^{{\rm Tr }_{\mathbb{F}_{2^m}/\mathbb{F}_2}(x+\overline{x})})(1-(-1)^{{\rm Tr}_{\mathbb{F}_{2^m}/\mathbb{F}_2}(x\overline{x})})\\
  &=&\frac{1}{4}(L_1-L_2-L_3+L_4),
\end{eqnarray*}
where
\begin{eqnarray*}
  &&L_1=\sum_{x\in \mathbb{F}_q^*}(-1)^{{\rm Tr }_{\mathbb{F}_{2^m}/\mathbb{F}_2}(ax+\overline{ax})}, L_2=\sum_{x\in \mathbb{F}_q^*}(-1)^{{\rm Tr }_{\mathbb{F}_{2^m}/\mathbb{F}_2}((a+1)x+(\overline{a}+1)\overline{x})} \\
 &&L_3=\sum_{x\in \mathbb{F}_q^*}(-1)^{{\rm Tr }_{\mathbb{F}_{2^m}/\mathbb{F}_2}(ax+\overline{a}\overline{x}+x\overline{x})}, L_4=\sum_{x\in \mathbb{F}_q^*}(-1)^{{\rm Tr }_{\mathbb{F}_{2^m}/\mathbb{F}_2}((a+1)x+(\overline{a}+1)\overline{x}+x\overline{x})}.
\end{eqnarray*}
It is easily seen that $\chi_1(D)=-\chi_0(D)$. Thus, in the next sequel, we assume that $a\neq 0,1$.
Now
\begin{equation*}
    L_1=\sum_{y\in \mathbb{F}_{2^m}^*, z\in \mathfrak{S}}(-1)^{{\rm Tr }_{\mathbb{F}_{2^m}/\mathbb{F}_2}(y(az+\overline{a}z^{-1}))}=1\cdot(2^m-1)+2^m\cdot(-1)=-1,
\end{equation*}
and if $a\neq 1$, then
\begin{equation*}
    L_2=\sum_{y\in \mathbb{F}_{2^m}^*, z\in \mathfrak{S}}(-1)^{{\rm Tr }_{\mathbb{F}_{2^m}/\mathbb{F}_2}(y((a+1)z+(\overline{a}+1)z^{-1}))}=1\cdot(2^m-1)+2^m\cdot(-1)=-1,
    \end{equation*}
For the next sum,
  \begin{equation*}
    L_3=\sum_{y\in \mathbb{F}_{2^m}^*, z\in \mathfrak{S}}(-1)^{{\rm Tr }_{\mathbb{F}_{2^m}/\mathbb{F}_2}(y(az+\overline{a}z^{-1}+1))}.
    \end{equation*}
    If $az+\overline{a}z^{-1}+1=0$ for some $z $ in an extension of $\mathbb{F}_q$, then putting $z=\frac{1}{a}u$, we have $u^2+u+a\overline{a}=0$. Since $a\overline{a}\in \mathbb{F}_{2^m}$, ${\rm Tr}(a\overline{a})=0$, the equation $X^2+X+a\overline{a}=0$ has two solutions in $\mathbb{F}_q$. If $u$ is a solution to the equation  $X^2+X+a\overline{a}=0$, then the other root of the equation is $\overline{u}\neq u$ if ${\rm Tr}_1^m(a\overline{a})=1$. Thus, $u\overline{u}=a\overline{a}, u+\overline{u}=1$. Hence we have $z\overline{z}=\frac{u\overline{u}}{a\overline{a}}=1$, and thus $z\in \mathfrak{S}$. So that $az+\overline{a}z^{-1}+1=0$ has exactly two distinct solutions $z\in \mathfrak{S}$. It is easy to see that $az+\overline{a}z^{-1}+1=0$ has no root in $\mathfrak{S}$ if ${\rm Tr}_1^m(a\overline{a})=0$. Therefore,
    \begin{equation}\label{f-28}
        L_3=\left\{\begin{array}{cc}
                    2^m-1, & \mbox{ if ${\rm Tr}_1^m(a\overline{a})=1$}, \\
                     (-1)(2^m+1), & \mbox{ if ${\rm Tr}_1^m(a\overline{a})=0$}.
                   \end{array}
        \right.
    \end{equation}

Finally, for $L_4$, using a similar method as before, we have
\begin{equation}\label{f-29}
L_4= \left\{\begin{array}{cc}
                    2^m-1, & \mbox{ if ${\rm Tr}_1^m((a+1)(\overline{a}+1))=1$}, \\
                     (-1)(2^m+1), & \mbox{ if ${\rm Tr}_1^m((a+1)(\overline{a}+1))=0$}.
                   \end{array}
        \right.
\end{equation}

Thus, we know that

(1) when $m$ is even, then
\begin{equation*}
   \chi_a(D)=\left\{\begin{array}{ll}
                     -2^{m-1}, & \mbox{ if ${\rm Tr}_1^n(a)=1$ and ${\rm Tr}_1^m(a\overline{a})=1$,} \\
                     2^{m-1},& \mbox{ if ${\rm Tr}_1^n(a)=1$ and ${\rm Tr}_1^m(a\overline{a})=0$,}\\
                     0,& \mbox{ in all other cases}.
                   \end{array}
   \right.
\end{equation*}

(2) when $m$ is odd, then
\begin{equation*}
   \chi_a(D)=\left\{\begin{array}{ll}
                     -2^{m-1}, & \mbox{ if ${\rm Tr}_1^n(a)=0$ and ${\rm Tr}_1^m(a\overline{a})=1$,} \\
                     2^{m-1},& \mbox{ if ${\rm Tr}_1^n(a)=0$ and ${\rm Tr}_1^m(a\overline{a})=0$,}\\
                     0,& \mbox{ in all other cases}.
                   \end{array}
   \right.
\end{equation*}

In summary, we have the following

{\thm \label{thm-4.5} Let $q=2^n$, $n=2m$ and let $\mathbb{F}_q$ be the finite field with $q$ elements. Define
\begin{equation*}
    D=\{x\in \mathbb{F}_q^*|{\rm Tr}_{\mathbb{F}_{2^m}/\mathbb{F}_2}(x+\overline{x})={\rm Tr}_{\mathbb{F}_{2^m}/\mathbb{F}_2}(x\overline{x})=1\}.
\end{equation*}
Then the Cayley graph $\mathfrak{C}(\mathbb{F}_q, D)$ is a connected bipartite graph which has only five eigenvalues, the eigenvalues of it are: $\pm 2^{n-2}$ if $m$ is even, and $\pm(2^{n-2}+2^{m-1})$ if $m$ is odd, and zero and $\pm 2^{m-1}$. As a consequence, $\mathfrak{C}(\mathbb{F}_q, D)$ is a Ramanujan Graph.}

Below we present an example to illustrate that one can obtain some nice graphs using Theorem \ref{thm-4.5}.

\noindent{\bf Example 4.6} Let $m=2$ and $q=2^n$, $n=2m$ and let $\mathbb{F}_q$ be the finite field with $q$ elements. Let $\gamma$ be a primitive element of $\mathbb{F}_q$ and $\mathfrak{C}(\mathbb{F}_q, D)$ the Cayley graph defined in Theorem \ref{thm-4.5}.
Then the neighbors of each vertex are listed in the following Table 4:
\begin{center}

\mbox{Table 4. Neighbors of the vertices}.

   \begin{tabular}{|c|c||c|c|}
     \hline
     \mbox{Vertex($\gamma^i$)}&\mbox{neighbors($\gamma^i$)}& \mbox{Vertex($\gamma^i$)}&\mbox{neighbors($\gamma^i$)} \\
     \hline
     \mbox{-$\infty$}&11, 13, 14, 7 & \mbox{0}&3, 6, 9, 12 \\
     \hline
     \mbox{1}&6, 7, 12, 14& \mbox{2}&9, 12, 13, 14 \\
      \hline
     \mbox{3}&0, 4, 5, 8& \mbox{4}&3, 9, 11, 13 \\
     \hline
      \mbox{5}&3, 7, 12, 13& \mbox{6}& 0, 1, 8, 10 \\
      \hline
      \mbox{7}&$-\infty$, 1, 5, 8& \mbox{8}& 3,  6, 7, 11 \\
      \hline
      \mbox{9}&0, 2, 4, 10& \mbox{10}& 6, 9, 11, 14 \\
      \hline
      \mbox{11}&$-\infty$, 4, 8, 10 & \mbox{12}&0, 1, 2, 5 \\
      \hline
      \mbox{13}&$-\infty$, 2, 4, 5 & \mbox{14}&$-\infty$, 1, 2, 10 \\
      \hline
   \end{tabular}
    \end{center}
   The spectral of $\mathfrak{C}(\mathbb{F}_q, D)$ are $\pm 4$ with multiplicity $1$, and $\pm 2$ with multiplicity $4$ and $0$ with multiplicity $6$, the diameter of $\mathfrak{C}(\mathbb{F}_q, D)$ is $4$. Moreover, $\mathfrak{C}(\mathbb{F}_q, D)$ is a distance transitive graph, the intersection array of it is $[4,3,2,1,1,2,3,4]$. The automorphism group of $\mathfrak{C}(\mathbb{F}_q, D)$ has order $384$.
\section*{Conclusion}
In this paper, we investigate the constructing for expander graphs, especially for Ramanujan graphs. Some infinite families of Ramanujan graphs are presented. It should be noted that the methods applied in this paper can be extend to more broader cases, for example, The subsets $C_0$ and $C_1$ in Theorem \ref{thm-33} may be any subgroup (delete the element $0$) in related groups, the resulting graphs may have some good properties, and the subset $D$ in Theorem \ref{thm-4.4} or in Theorem \ref{thm-4.5} may be some algebraic varieties. If one can determine the character sum (a kind of Weil sums) of the set, then he can determine some parameters of the graph. Some further combinatorial and cryptographic properties of the constructed graphs in this paper will be discussed in another paper.

\section*{Acknowledgement}
For the graph in Theorem \ref{thm-4.5}, we have a discussion with Prof. Qing Xiang, Prof. C. Carlet, Prof. Cunsheng Ding et al. We should express our grateful thankfulness to all of them for valuable comments. The work of this paper is supported by the NUAA Fundamental Research Funds, No. 2013202.


\end{document}